\documentclass[reqno]{amsart}
\usepackage[T1]{fontenc}
\usepackage[latin1]{inputenc}

\makeatletter

\providecommand{\LyX}{L\kern-.1667em\lower.25em\hbox{Y}\kern-.125emX\@}

    \theoremstyle{plain}    
    \newtheorem{cor}{Corollary}[section] 

 \theoremstyle{plain}    
 \newtheorem{thm}{Theorem} 
 \theoremstyle{plain}    
 \newtheorem{lem}[cor]{Lemma} 
 \theoremstyle{plain}    
 \newtheorem{prop}[cor]{Proposition} 
 \theoremstyle{definition}
 \newtheorem{defn}{Definition}
 \theoremstyle{remark}
 \newtheorem*{rem*}{Remark}

\newcommand{\ZZ}{\mathbb {Z}}
\newcommand{\NN}{\mathbb {N}}
\newcommand{\TT}{\mathbb {T}}
\newcommand{\RR}{\mathbb {R}}
\newcommand{\CC}{\mathbb {C}}

\newcommand{\one}{\mathbf {1}}

\newcommand{\spec}{\mathrm{spec}\,}

\newcommand{\sign}{\mathrm{sign}\,}
\newcommand{\dist}{\mathrm{dist}\,}

\newcommand{\half}{{\textstyle \frac{1}{2}}}
\newcommand{\third}{{\textstyle \frac{1}{3}}}
\newcommand{\twothirds}{{\textstyle \frac{2}{3}}}
\newcommand{\quarter}{{\textstyle \frac{1}{4}}}
\newcommand{\threequarters}{{\textstyle \frac{3}{4}}}

\newcommand{\brk}{\\ }

\makeatother

\begin{document}

\title{Menshov Representation Spectra}

\author{Gady Kozma and Alexander Olevski\v\i}

\address{School of Mathematics, Tel Aviv University, Israel}

\email{gady@math.tau.ac.il, olevskii@math.tau.ac.il}

\maketitle

\section{\label{section_intro}Introduction}

\subsection{{}}

The problem of representing a function \( f \) on the circle \( \TT  \) by
a trigonometric series
\begin{equation}
\label{f_is_sum}
f(t)=\sum _{k\in \ZZ }c_{k}e^{ikt}
\end{equation}
has a long history. Most of the attention has been paid to Fourier expansions,
where the coefficients are derived from \( f \) using an integration process.
Certainly this approach requires \'{a} priori the integrability of \( f \)
and even under this assumption it does not always ensure satisfactory convergence
properties. 

On the other hand, D. E. Menshov \cite{11} discovered that every \( f\in L^{0}(\TT ) \),
that is any measurable function from \( \TT  \) to \( \CC  \), with no additional
requirements, can be represented as a sum of a series (\ref{f_is_sum}) converging
almost everywhere (a.e.). The coefficients here are obtained by a special construction
which has certain freedoms (for example, one can avoid using low frequency oscillations). 

This result (usually called ``Menshov representation theorem'') gave rise
to many further investigations, aimed to find in which directions it can be
extended and improved and in which it cannot. Many interesting results on this
problem were obtained by Menshov himself \cite{12}, N. Bary \cite{ba} (chap.
XV), Talalyan \cite{tal}, Arutyunyan \cite{ayan}, Kashin \cite{13}, Konyagin
\cite{14}, K\"{o}rner \cite{ko} and others. The reader may find a comprehensive
bibliography (up to '92) in the survey \cite{survey}. We add also our recent
paper \cite{Kozol}.

\subsection{{}}

One aspect of the theory is: how rich must the family of exponentials be to
ensure the possibility of representation. We introduce the following

\begin{defn}
A set \( \Lambda  \) of integers, \( \Lambda =\{\lambda (n)\, ;\, \ldots <\lambda (-1)<\lambda (0)<\lambda (1)<\ldots \} \)
is a \emph{Menshov spectrum} if every \( f\in L^{0}(\TT ) \) admits a representation
\begin{equation}
\label{f_is_sum_on_lambda}
f(t)=\sum _{k\in \Lambda }c_{k}e^{ikt}\equiv \lim _{N\rightarrow \infty }\sum _{k\in \Lambda ,|k|\leq N}c_{k}e^{ikt}
\end{equation}
converging a.e.
\end{defn}
It is well known that a Menshov spectrum might be quite sparse. For example,
it may have density zero. Arutyunyan proved \cite{ayan} that any \( \Lambda  \)
symmetric with respect to zero and containing arbitrarily long segments of integers
is a Menshov spectrum. 

The first result in our paper shows that a Menshov spectra might even be lacunary,
that is satisfy the condition
\[
\lambda (n+1)-\lambda (n)\rightarrow \infty \quad (|n|\rightarrow \infty )\]
 and the sizes of the gaps may grow quite fast.

\begin{thm}
\label{Theorem_sparse}Given a positive sequence \( \epsilon (n)=o(1) \) as
\( n\rightarrow \infty  \), one can define a symmetric Menshov spectrum \( \Lambda  \)
such that
\begin{equation}
\label{lambda_near_exponent}
\frac{\lambda (n+1)}{\lambda (n)}>1+\epsilon (n)\quad n=1,2,\ldots 
\end{equation}
 
\end{thm}
This result is sharp: the \( \epsilon (n) \) on the right hand side cannot
be replaced by a positive constant, which would mean that \( \Lambda  \) is
lacunary in the sense of Hadamard. Indeed, it is well known that if for such
a \( \Lambda  \) the series (\ref{f_is_sum_on_lambda}) converges then \( f \)
must satisfy some special properties, see for example \cite{ba} (chap. XI)
or \cite{17} (chap. V).

\subsection{{}}

These restrictions on the sizes of the gaps are far from giving a complete picture.
The arithmetics of the set play a crucial role as well. It is obvious, for example,
that \( \Lambda =2\ZZ  \) is not a Menshov spectrum since the sum (\ref{f_is_sum_on_lambda})
has in this case a period \( <2\pi  \). Somewhat more delicate arguments show
that the set \( \Lambda =\{\pm k^{2}\} \) is not a Menshov spectrum, even after
any bounded perturbation. On the other hand, we prove the following: 

\begin{thm}
\label{theorem_squares}For any sequence \( w(k)\rightarrow \infty  \) as \( k\rightarrow \infty  \)
one can construct a symmetric Menshov spectrum \( \Lambda  \),
\begin{equation}
\label{lambda_is_squares}
\Lambda =\left\{ \pm k^{2}+o(w(|k|))\right\} \quad .
\end{equation}

\end{thm}
We do not know whether such a result for cubes is true.

\subsection{{}}

Consider now the case when \emph{only positive} frequencies are involved. For
obvious reasons this situation can be referred to as ``analytic'': having
a decomposition (\ref{f_is_sum_on_lambda}) with \( \Lambda =\ZZ ^{+} \) one
can consider an analytic function in the unit disc 
\[
F(z)=\sum _{k>0}c_{k}z^{k}\quad .\]
 At any point \( t\in \TT  \) where the sum (\ref{f_is_sum_on_lambda}) converges,
according to the classical Abel theorem, \( F(z)\rightarrow f(t) \) when \( z \)
approaches the point \( e^{it} \) non-tangentially. The Lusin-Privalov uniqueness
theorem, see for example \cite{koosis}, says that ``non-tangential boundary
values'' of an analytic function can not be an arbitrary function, specifically,
if it vanishes on a set of positive measure then it is identically zero. It
is interesting that for \emph{radial} limits a.e. there are no such limitations:
Kahane and Katznelson \cite{kk} proved that any measurable function \( f \)
might be represented in such a form. 

So the argument above leads us to the conclusion that \( \ZZ ^{+} \) \emph{is
not a Menshov spectrum}. 

It turns out, however, that if we relax a bit the notion of convergence, by
replacing pointwise convergence with convergence in \( L^{0} \) metric, we
might again get an unrestricted representation theorem

\begin{thm}
\label{Theorem_L0}Every \( f\in L^{0}(\TT ) \) can be represented as a sum
\begin{equation}
\label{f_is_sum_measure}
f(t)=\sum _{k>0}c_{k}e^{ikt}
\end{equation}
converging in measure.
\end{thm}
It should be mentioned that Menshov \cite{12} showed that using convergence
in measure one may get a representation (\ref{f_is_sum}) even for functions
\( f:\TT \rightarrow \RR \cup \{\pm \infty \} \). Theorem \ref{Theorem_L0}
is also valid for such functions, so it gives a direct improvement of Menshov's
theorem by replacing the whole spectrum \( \ZZ  \) by \( \ZZ ^{+} \).

Finally, we prove that relatively ``small'' parts of \( \ZZ ^{+} \) are already
sufficient for representation in measure. In particular we construct ``almost
Hadamarian'' and ``almost squares'' Menshov spectra in measure.

Theorems \ref{Theorem_sparse} and \ref{theorem_squares} are proved below in
section \ref{section_lacunary}. Theorem \ref{Theorem_L0} (in a stronger form)
is proved in section \ref{section_analytic}, and the last mentioned results
in section \ref{section_analytic_lacunary}.

The main results of this paper were stated without proofs in our note \cite{kozol_short}.

\section{Preliminaries}

\subsection{Notations}

We denote by \( \TT  \) the circle group, \( \TT =\RR /2\pi \ZZ  \), and we
identify it in the standard way with the segment \( [-\pi ,\pi ] \). 

\( \mathbf{m} \) is the normalized Lebesgue measure on \( \TT  \). 

\( \ZZ ^{+} \) or \( \NN  \) is the set \( \{1,2,...\} \).

\( C \) is used to denote absolute positive constants, possibly different.

For a set \( A\subset \TT  \) we define by \( \one _{A} \) the function 
\[
\one _{A}(t)=\left\{ \begin{array}{ll}
1 & t\in A\\
0 & t\notin A
\end{array}\right. \]

The \( L_{0} \) ``norm'' of any measurable function \( f \) from \( \TT  \)
to \( \CC  \) is defined to be 
\begin{equation}
\label{L0_metric}
\left\Vert f\right\Vert _{0}:=\inf \left\{ \epsilon \, :\, \mathbf{m}\{t\, :\, |f(t)|>\epsilon \}<\epsilon \right\} \quad .
\end{equation}
 This norm is sub-additive but not homogeneous. \( L_{0}(\TT ) \) is the space
of all measurable functions on the segment \( \TT  \) endowed with the norm
\( ||\cdot ||_{0} \) and is a topological linear space. Convergence in \( ||\cdot ||_{0} \)
is equivalent to convergence in measure.

For \( f\in L^{1}(\TT ) \) we denote by \( \widehat{f}(n) \) the Fourier coefficients
of \( f \) and write 
\[
f(t)\sim \sum \widehat{f}(k)e^{ikt}\quad .\]
 The spectrum of \( f \) (denoted by \( \spec f \)) is the support of \( \widehat{f} \),
i.e. the set of integers where \( \widehat{f} \) is non zero.

The triangle function of width \( 2\epsilon  \), denoted by \( \tau _{\epsilon } \)
is defined on \( \left[ -\pi ,\pi \right]  \) by
\begin{equation}
\label{def_tau_delta}
\tau _{\epsilon }=\left\{ \begin{array}{cc}
1-\frac{|x|}{\epsilon } & |x|<\epsilon \\
0 & \mathrm{otherwise}
\end{array}\right. \quad .
\end{equation}

It is well known (and easy to calculate) that \( \widehat{\tau _{\epsilon }}(n)>0 \)
for any \( n\in \ZZ  \) and therefore 
\begin{equation}
\label{tau_A1_1}
||\widehat{\tau _{\epsilon }}||_{1}=1
\end{equation}
 (as usual, \( ||\cdot ||_{p} \), \( 1\leq p\leq \infty  \), for a function
or a sequence denotes the \( L^{p}(\TT ) \) or \( l_{p} \) norms respectively).
A trigonometric polynomial is a finite sum
\[
P(t)=\sum c_{k}e^{ikt}\]
the degree of \( P \), denoted by \( \deg P \) is \( \max \{|k|\, :\, k\in \spec P\} \).
We say about two trigonometric polynomials that the spectrum of \( P \) \label{def_follows}follows
the spectrum of \( Q \) (or simply that \( P \) follows \( Q \)) if
\[
l\in \spec Q,\; k\in \spec P\Rightarrow |l|<|k|\quad .\]

For a trigonometric polynomial \( P=\sum c_{k}e^{ikt} \) we denote the partial
sums by 
\[
S_{n}(P;t):=\sum _{k=-n}^{n}c_{k}e^{ikt}\quad .\]
and their maximum by 
\[
S^{*}(P;t):=\sup _{n\geq 0}\left| S_{n}(P;t)\right| ,\]
Occasionally, we shall also need estimates of non-symmetric partial sums, so
define
\[
S_{n,m}(P;t):=\sum _{k=m}^{n}c_{k}e^{ikx}\]
\[
S^{**}:=\sup _{n,m\in \ZZ }\left| S_{n,m}\right| \quad .\]
Finally, we shall often use the following obvious, and somewhat rough estimates:
\[
||S^{*}(P)||_{\infty },\, ||S^{**}(P)||_{\infty }\leq ||\widehat{P}||_{1}\]

\subsection{\label{ssec_cross}Special Products}

Given \( f:\TT \rightarrow \CC  \) and \( r\in \NN  \) we denote by \( f_{[r]} \)
the ``contracted'' function
\[
f_{[r]}(t)=f(rt)\quad t\in \TT .\]
We shall systematically use through this paper a ``special product'' of the
form 
\begin{equation}
\label{cross}
H=Q_{[r]}\cdot P\quad .
\end{equation}
 where \( P \) and \( Q \) are trigonometric polynomials, \( r>2\deg P \)
and \( \widehat{Q}(0)=0 \). The spectrum of \( H \) has a block structure,
with each block containing a translated copy of \( \spec P \). So, if \( n>0 \):
\[
n=sr+l,\quad -\half r\leq l<\half r\]
we have:
\begin{equation}
\label{Sn_of_PcrossQ}
\sum _{j=0}^{n}\widehat{H}(j)e^{ijt}=P(t)\cdot \sum _{k=1}^{s-1}\widehat{Q}(k)e^{ikrt}+\widehat{Q}(s)e^{isrt}\cdot \sum _{j=-\deg P}^{l}\widehat{P}(j)e^{ijt}\quad .
\end{equation}
If \( n<0 \) an analogous equality holds for \( \sum _{j=n}^{0} \). 

The second sum on the right hand side we usually bound by \( ||\widehat{Q}||_{\infty }\cdot ||\widehat{P}||_{1} \).

Sometimes, as in the proof of lemma \ref{Lemma_Korner} below, we use the trivial
estimate for \( H \):
\begin{equation}
\label{Sstartstar_H_lt_PQ}
||S^{**}(H)||_{\infty }\leq ||\widehat{P}||_{1}\cdot ||\widehat{Q}||_{1}
\end{equation}

\subsection{Separation of spectra}

Our starting point is a fundamental principle, used by Menshov in his ``correction''
and ``representation'' theorems. In a few words, it goes as follows: given
\( f \) we decompose it into ``elementary pieces'' and then by correcting
in a small measure we can localize the Fourier transform of each piece in a
``personal'' part of the spectrum in such a way that the supports of different
pieces are essentially disjoint. It allows to avoid ``unpleasant'' resonance
with the Dirichlet kernel and control the size of the Fourier partial sums.
Menshov did it using a delicate and difficult technique based on his so-called
``Dirichlet factor lemma'' (see \cite{ba}, chap. V). A new version of Menshov's
approach (using the special products above) was outlined in \cite{ol}, sect.
1, in the framework of the correction theorem. A similar technique for representation
results was developed by Körner, see \cite{ko}. The following lemma may serve
as a good illustration of the modern version of Menshov's approach:

\begin{lem}
\label{Lemma_Korner}Given any \( \epsilon >0 \), \( \delta >0 \), there exists
a trigonometric polynomial \( Q \) such that
\end{lem}
\begin{enumerate}
\item \label{req_Qhat_lem_korner}\( \widehat{Q}(0)=0 \), \( ||\widehat{Q}||_{\infty }<\delta  \);
\item \label{req_Q1_lem_korner}\( \mathbf{m}\left\{ t\, :\, |Q(t)-1|>\delta \right\} <\epsilon  \);
\item \label{req_Sq_lem_korner}\( ||S^{**}(Q)||_{\infty }<C\epsilon ^{-1} \).
\end{enumerate}
Such a statement, with a little weaker form of \ref{req_Sq_lem_korner} can
be found in \cite{ko}, lemma 1; but we need the present form, especially for
the proof of theorem \ref{theorem_squares}. The proof below follows essentially
\cite{ol}, sect. 1.

\begin{proof}
Fix a large integer \( K \). Let \( f:=\tau _{2\pi /K} \), \( g:=1-\frac{4}{\epsilon }\tau _{\epsilon /4} \).
Approximating them by Fourier partial sums we get polynomials \( F \) and \( G \)
satisfying
\[
||\widehat{F-f}||_{1}<c,\quad ||\widehat{G-g}||_{1}<c,\quad \widehat{G}(0)=0\]
where \( c \) is small. Let \( T \) be the translation by \( \frac{2\pi }{K} \).
Set 
\begin{eqnarray}
Q_{s} & := & T^{s}(F)\cdot G_{[N_{s}]};\label{V2_Ns} \\
Q & := & \sum _{s=1}^{K}Q_{s}\quad ,\nonumber 
\end{eqnarray}
where \( N(1)>2\deg F \) and \( N(s) \) grow sufficiently fast to ensure the
separation of spectra. It follows:
\[
||\widehat{Q}||_{\infty }=\max _{s}||\widehat{Q_{s}}||_{\infty }=||\widehat{F}||_{\infty }||\widehat{G}||_{1}<\frac{2}{K}\cdot \frac{4}{\epsilon }\quad .\]
So, if \( K(\epsilon ,\delta ) \) is sufficiently large we get \ref{req_Qhat_lem_korner}.
Obviously 
\[
\sum _{s=1}^{K}T^{s}(f)\cdot g_{[N(s)]}=1\]
outside a set of measure \( <\epsilon  \), so if \( c \) is sufficiently small
we get \ref{req_Q1_lem_korner}. Now denote:
\[
\widetilde{G}(t)=\sum _{k>0}\widehat{G}(k)e^{ikt}\]
The separation of spectra implies a decomposition 
\[
\Sigma ^{(n)}:=\sum _{k=0}^{n}\widehat{Q}(k)e^{ikt}=\Sigma _{1}+\Sigma _{2}\]
where \( \Sigma _{1} \) contains all polynomials of the form \( T^{s}(F)\widetilde{G}_{[N(s)]} \)
with spectra belonging to \( [0,n] \) and the reminder \( \Sigma _{2} \) is
a segment of a polynomial of the same form. According to section \ref{ssec_cross}
we have
\begin{eqnarray*}
||\Sigma _{2}||_{\infty } & \leq  & ||\widehat{F}||_{1}\cdot ||\widehat{G}||_{1}\\
||\Sigma _{1}||_{\infty } & \leq  & \left\Vert \sum |T^{s}(F)|\right\Vert _{\infty }\cdot ||\widehat{G}||_{1}\quad ,
\end{eqnarray*}
so \( c<\frac{1}{K} \) implies \( ||\Sigma ^{(n)}||_{\infty }<C\epsilon ^{-1} \).
The same is true for \( n<0 \) and we get \ref{req_Sq_lem_korner}.
\end{proof}

\subsection{Analytic polynomials.}

An analytic polynomial is a trigonometric polynomial \( P \) with positive
spectrum, i.e.
\[
P(t)=\sum _{n>0}c_{n}e^{int}\quad .\]
Sections \ref{section_analytic} and \ref{section_analytic_lacunary} will contain
most of the discussion about these objects. For the time being, we wish to remind
the reader the generally known fact that analytic polynomials are dense in \( L^{0} \).
We formulate it in an equivalent way:

\begin{lem}
\label{Q_eps}For every \( \epsilon >0 \) there exists an analytical polynomial
\( R_{\epsilon } \) such that 
\begin{equation}
\label{Q_is_1}
\left\Vert R_{\epsilon }-1\right\Vert _{0}<\epsilon 
\end{equation}

\end{lem}
\begin{proof}
Let \( g(e^{it}) \) be a real (smooth) function with zero average and \( <-1/\epsilon  \)
for \( |t|>\epsilon /3 \). Denote by \( G \) the Poisson integral of \( g \);
by \( \tilde{G} \) the conjugate harmonic function and write
\[
F=\exp \left( G+i\tilde{G}\right) \]
Clearly, \( F(0)=1 \) and \( ||F(e^{it})||_{0}<\epsilon  \) so taking \( R_{\epsilon } \)
to be a Taylor approximation of \( 1-F \) we get the result.
\end{proof}

\section{\label{section_lacunary}Lacunary Spectra}

This section is devoted to the proof of theorems \ref{Theorem_sparse} and \ref{theorem_squares}.

\subsection{The blocks \protect\( B(s,a)\protect \). }

Our starting point to these two theorems, as well as to the theorems of section
\ref{section_analytic_lacunary}, are blocks \( B(s,a)\subset \ZZ  \), defined
for \( s \), \( a\in \NN  \) as follows:
\begin{eqnarray}
B_{1}(s,a) & := & \{-sa,-(s-1)a,\ldots ,-a,a,2a,\ldots ,sa\}\nonumber \\
B^{+}_{1}(s,a) & := & B_{1}(s,a)\cap \ZZ ^{+}\nonumber \\
B_{2}(s,a) & := & \bigcup _{|k|\leq s}k+B^{+}_{1}(s,(2s)^{k+s}a)\label{def_Bna2} \\
B(s,a) & := & B_{1}(s,(2s)^{2s+2}a)+\left( B_{2}(s,a)\cup -B_{2}(s,a)\right) \nonumber \label{how_can_I_remove_this} 
\end{eqnarray}
where we understand \( x+B \) and \( A+B \) in the standard sense: 
\begin{eqnarray*}
x+A & := & \left\{ x+a\, :\, a\in A\right\} \\
A+B & := & \{a+b\, :\, a\in A,\, b\in B\}\quad .
\end{eqnarray*}

The following can be checked directly from the definition:

\begin{prop}
\label{Bsa_linear}For any \( b\in B(s,a) \) there exists \( l=l(b) \), \( 0<|l|<C(s) \)
such that \( |b-la|<C(s) \). Further, \( b_{1}\neq b_{2}\Rightarrow l(b_{1})\neq l(b_{2}) \).
\end{prop}
This means that \( B(s,a) \) is a subset of a ``small'' perturbation of a
linear progression (if \( a\gg s \)). In particular, \( B(s,a) \) is quite
sparse: for \( b_{1}\neq b_{2} \) in \( B(s,a) \) 
\begin{equation}
\label{Bsa_sparse}
|b_{1}-b_{2}|>a-2C(s)\quad .
\end{equation}
Another point to notice is that \( B(s,a) \) has a rather large ``hole''
near zero. We might quantify this as 
\begin{equation}
\label{B_hole_at_0}
\dist (B,0)>(2s)^{2s+1}a\quad .
\end{equation}

\( B(s,a) \) is clearly symmetric with respect to zero. Actually, the component
\( -B_{2}(s,a) \) is here only to assure that. For the proof of lemma \ref{Lemma_P_in_Bna}
below, it would have been enough to define \( B(s,a) \) as \( B_{1}(s,(2s)^{2s+2}a)+B_{2}(s,a) \)

\subsection{Approximation with spectrum in \protect\( B(s,a)\protect \).}

\begin{lem}
\label{Lemma_P_in_Bna}For every \( \epsilon >0 \), \( \delta >0 \) and \( f\in L^{0} \)
there exists an \( S=S(f,\epsilon ,\delta ) \) with the following property.
Given \( s>S \) and \( a\in \NN  \) one may construct a trigonometric polynomial
\( P \) satisfying 
\begin{enumerate}
\item \label{ref_Pisf_lem_Bsa}\( \mathbf{m}\left\{ t\, :\, |P-f|>\delta \right\} <\epsilon  \);
\item \label{ref_specP_lem_Bsa}\( \spec P\subset B(s,a) \);
\item \label{ref_SP_lem_Bsa}\( \mathbf{m}\left\{ t\, :\, S^{**}(P;t)>C\epsilon ^{-1}(|f(t)|+\delta )\right\} <\epsilon  \).
\end{enumerate}
\end{lem}
\begin{proof}
We shall perform this approximation in three steps. 

\emph{Step 1:} Approximate \( f \) by a polynomial \( P_{1} \) satisfying
\begin{equation}
\label{f_minus_p1}
\mathbf{m}\left\{ t\, :\, |f-P_{1}|>\delta _{1}\right\} <\epsilon _{1}
\end{equation}
where
\[
\delta _{1}:=\third \delta ,\quad \epsilon _{1}:=\third \epsilon \quad .\]

\emph{Step 2:} Define 
\begin{equation}
\label{eps_2_del_2}
\epsilon _{2}:=\frac{\epsilon }{6\deg P_{1}+3},\quad \delta _{2}:=\frac{\delta }{3||\widehat{P_{1}}||_{1}}
\end{equation}
 And use lemma \ref{Q_eps} to get an analytic polynomial \( Q_{2} \) satisfying
\begin{equation}
\label{Q2_is_1}
\mathbf{m}\left\{ t\, :\, |Q_{2}-1|>\delta _{2}\right\} <\epsilon _{2}
\end{equation}

\emph{Step 3:} Use lemma \ref{Lemma_Korner} with 
\begin{equation}
\label{eps_3_del_3}
\epsilon _{3}:=\third \epsilon ,\quad \delta _{3}:=\frac{\delta }{6||\widehat{P_{1}}||_{1}||\widehat{Q_{2}}||_{1}}
\end{equation}
to get \( Q_{3} \) satisfying

\begin{eqnarray}
 &  & \widehat{Q_{3}}(0)=0;\label{Q3hat_zero_zero} \\
 &  & ||\widehat{Q_{3}}||_{\infty }\leq \delta _{3};\label{Q3_hat_is_small} \\
 &  & ||S^{**}(Q_{3})||_{\infty }\leq C\epsilon _{3}^{-1};\label{Sstar_of_Q3} \\
 &  & \mathbf{m}\left\{ t\, :\, |Q_{3}-1|>\delta _{3}\right\} <\epsilon _{3}\quad .\nonumber 
\end{eqnarray}

We have now gathered all the data we need to define \( S \):
\begin{equation}
\label{def_S_eps_del_f}
S(f,\epsilon ,\delta ):=\max \left\{ \deg P_{1},\deg Q_{2},\deg Q_{3}\right\} \quad .
\end{equation}
Let us now show how we can construct \( P \) supported on \( B(s,a) \). Fix
some \( s>S \) and any integer \( a \). Clearly (remember that \( Q_{2} \)
is analytic) 
\[
\spec (Q_{2})_{[a]}\subset B^{+}_{1}(s,a)\]
or
\begin{equation}
\label{supp_exp_R}
\spec \left( e^{ikt}\cdot (Q_{2})_{[a]}\right) \subset k+B^{+}_{1}(s,a).
\end{equation}
So define 
\begin{eqnarray}
P_{2} & := & \sum _{|k|\leq \deg P_{1}}\widehat{P_{1}}(k)e^{ikt}\cdot (Q_{2})_{[p(k)]}\label{def_P2} \\
p(k) & := & a(2s)^{k+s}\label{def_pk} 
\end{eqnarray}
 On one hand, (\ref{supp_exp_R}) and (\ref{def_Bna2}) imply: \( \spec P_{2}\subset B_{2}(s,a) \).
On the other hand,
\[
|P_{2}-P_{1}|=\left| \sum \widehat{P_{1}}(k)e^{ikt}(1-(Q_{2})_{[p(k)]})\right| \]
so (\ref{eps_2_del_2}) and (\ref{Q2_is_1}) gives
\begin{eqnarray}
\lefteqn {\mathbf{m}\left\{ t\, :\, |P_{2}-P_{1}|>\third \delta \right\} \leq } &  & \nonumber \\
 & \qquad  & \leq \sum _{|k|\leq \deg P_{1}}\mathbf{m}\left\{ t\, :\, \left| \widehat{P_{1}}(k)(1-(Q_{2})_{[p(k)]})\right| >|\widehat{P_{1}}(k)|\delta _{2}\right\} \label{P1_minus_P2} \\
 &  & <(2\deg P_{1}+1)\epsilon _{2}=\third \epsilon \nonumber 
\end{eqnarray}
 and thus
\begin{equation}
\label{f_minus_P2}
\mathbf{m}\left\{ t\, :\, |P_{2}-f|>\twothirds \delta \right\} <\twothirds \epsilon \quad .
\end{equation}
 Finally, define
\[
P\equiv P_{3}=(Q_{3})_{[p(s+2)]}P_{2}\quad .\]
From (\ref{Q3hat_zero_zero}) we get that \( \spec (Q_{3})_{[p(s+2)]}\subset B_{1}(s,p(s+2)) \),
which implies \ref{ref_specP_lem_Bsa}. First, an estimate of \( P-P_{2} \).
Obviously
\[
\mathbf{m}\left\{ x\, :\, |P-P_{2}|>\third \delta \right\} \leq \mathbf{m}\left\{ t\, :\, |1-Q_{3}|>\frac{\delta }{3||P_{2}||_{\infty }}\right\} \quad .\]
The definition of \( P_{2} \), (\ref{def_P2}), gives the estimate 
\[
||P_{2}||_{\infty }\leq ||\widehat{P_{1}}||_{1}||\widehat{Q_{2}}||_{1}\]
so, using (\ref{eps_3_del_3}),
\begin{equation}
\label{P3_minus_P2}
\mathbf{m}\left\{ t\, :\, |P-P_{2}|>\third \delta \right\} \leq \mathbf{m}\left\{ t\, :\, |1-Q_{3}|>\delta _{3}\right\} <\third \epsilon \quad .
\end{equation}
Summing (\ref{f_minus_P2}) and (\ref{P3_minus_P2}) gives
\[
\mathbf{m}\left\{ t\, :\, |f-P|>\delta \right\} <\epsilon \quad .\]
 To estimate \( S^{**}(P) \) note that, using (\ref{def_pk}), (\ref{def_P2})
and (\ref{def_S_eps_del_f}),
\begin{equation}
\label{p_s_plus_2_large}
p(s+2)=(2s)^{2s+2}a>2\deg P_{2},
\end{equation}
so we may use (\ref{Sn_of_PcrossQ}) and get
\[
S^{**}(P;t)\leq ||S^{**}(Q_{3})||_{\infty }|P_{2}(t)|+2||\widehat{Q_{3}}||_{\infty }||\widehat{P_{2}}||_{1}\quad .\]
As before, estimate \( ||\widehat{P_{2}}||_{1} \) directly from (\ref{def_P2})
and get 
\[
||\widehat{P_{2}}||_{1}\leq ||\widehat{P_{1}}||_{1}||\widehat{Q_{2}}||_{1}\]
which, together with (\ref{Sstar_of_Q3}), (\ref{Q3_hat_is_small}) and (\ref{eps_3_del_3})
gives
\[
S^{**}(P;t)\leq C\epsilon ^{-1}|P_{2}(t)|+\third \delta \quad .\]
and remembering (\ref{f_minus_P2}), the lemma is proved.
\end{proof}

\subsection{Almost Hadamarian lacunarity}

The proof of theorem \ref{Theorem_sparse} contains two components: showing
that every set \( \Lambda  \) containing \( B(s,a) \) for arbitrarily large
\( s \) is a Menshov spectrum and that such a set conforming to the requirements
of the theorem can be build. Let us start with the first component.

\begin{lem}
\label{Bna_is_menshov}If \( \Lambda \subset \ZZ  \) satisfies \( B(s_{k},a_{k})\subset \Lambda  \)
for a sequence \( s_{k}\rightarrow \infty  \) then it is a Menshov spectrum.
\end{lem}
\begin{proof}
Let \( f \) be in \( L^{0} \). We shall construct by induction a sequence
of polynomials \( P_{n} \) satisfying
\begin{enumerate}
\item \label{ref_fissum_lem_Bsa2}\( \mathbf{m}\left\{ x\, :\, |f-\sum _{k=1}^{n}P_{k}|>2^{-2n}\right\} <2^{-n} \);
\item \label{ref_specP_lem_Bsa2}\( \spec P_{n}\subset \Lambda  \);
\item \label{ref_Pfollow_lem_Bsa2}\( P_{n} \) follows \( P_{n-1} \);
\item \label{ref_SP_lem_Bsa2}\( ||S^{*}(P_{n})||_{0}\leq C2^{-n} \).
\end{enumerate}
This finishes the proof as it is clear that \( \sum _{n=1}^{\infty }P_{n} \)
yields a trigonometric series converging to \( f \) almost everywhere. 

Assume \( P_{0},\ldots ,P_{n-1} \) has been constructed. Define 
\begin{eqnarray*}
 &  & F_{n}:=f-\sum _{k=1}^{n-1}P_{k}\\
 &  & \epsilon _{n}:=2^{-n}\\
 &  & \delta _{n}:=2^{-2n}\quad .
\end{eqnarray*}
Lemma \ref{Lemma_P_in_Bna} gives \( S=S(F_{n},\epsilon _{n},\delta _{n}) \)
and we may pick some \( s>\max \left\{ S,\deg P_{n-1}\right\}  \) and some
\( a \) for which \( B(s,a)\subset \Lambda  \). Thus there exist \( P_{n} \)
satisfying the requirements of lemma \ref{Lemma_P_in_Bna}, namely 
\begin{eqnarray}
 &  & \mathbf{m}\left\{ t\, :\, |P_{n}-F_{n}|>\delta _{n}\right\} <\epsilon _{n}\label{Pn_minus_Fn} \\
 &  & \mathbf{m}\left\{ t\, :\, S^{*}(P_{n};t)>C\epsilon _{n}^{-1}(|F_{n}(t)|+\delta _{n})\right\} <\epsilon _{n}\label{Sstar_of_Pn} \\
 &  & \spec P\subset B(s,a)\subset \Lambda \quad .\nonumber 
\end{eqnarray}

Now, \( s>\deg P_{n-1} \) implies requirement \ref{ref_Pfollow_lem_Bsa2} of
the lemma (remember (\ref{B_hole_at_0})). Requirement \ref{ref_fissum_lem_Bsa2}
is an immediate consequence of (\ref{Pn_minus_Fn}) since 
\[
f-\sum _{k=1}^{n}P_{k}=F_{n}-P_{n}\]
requirement \ref{ref_SP_lem_Bsa2} of the lemma is a combination of (\ref{Sstar_of_Pn})
and our induction assumption (requirement \ref{ref_fissum_lem_Bsa2}) for \( P_{n-1} \):
\[
\mathbf{m}\left\{ x\, :\, |F_{n}|>2^{-2n+2}\right\} <2^{-n+1}\]
and the lemma is proved.
\end{proof}
To finish the proof of theorem \ref{Theorem_sparse} it is enough to find a
symmetric set \( \Lambda  \) conforming to (\ref{lambda_near_exponent}) and
containing a sequence of \( B(s_{k},a_{k}) \). It can be done by a simple induction.
We may assume that \( \{\epsilon _{n}\} \) is monotone. At each step we either
add to \( \Lambda  \) and entire block \( B \) or define two (symmetric) members
of \( \Lambda  \). More precisely, let
\[
\{\lambda (j)\}_{|j|\leq m}\]
 be already defined. 

(a) If there exists an \( s \) such that 
\begin{equation}
\label{eps_suffic_small}
\epsilon (m)<\frac{1}{2C(s)}\quad ,
\end{equation}
 where \( C(s) \) is the constant from proposition \ref{Bsa_linear}, we take
the maximal \( s \) satisfying this inequality and include all elements of
the block \( B(s,a) \) as new members of \( \Lambda  \). \( a \) here is
chosen sufficiently large to ensure (\ref{lambda_near_exponent}) for \( n=m \)
and simultaneously 
\[
\frac{a-2C(s)}{aC(s)}>\frac{1}{2C(s)}\]
which implies (\ref{lambda_near_exponent}) for all 
\[
m<n<m+\half |B(s,a)|\]
due to \( \max B(s,a)<C(s)a \) and (\ref{Bsa_sparse}).

(b) If such an \( s \) does not exist, we just define 
\[
\lambda (m+1)>\lambda (m)(1+\epsilon (m))\quad .\]
Obviously, the condition (\ref{eps_suffic_small}) will be fulfilled infinitely
many times, and the corresponding values of \( s \) will tend to infinity;
this finishes the proof.\qed

As discussed in the introduction, this result is sharp.

\subsection{Squares.}

Lemma \ref{Bna_is_menshov} allows one to get a Menshov spectrum satisfying
\[
\lambda (n)=\pm n^{2}+O(n)\quad .\]
To obtain such a result with \emph{arbitrarily slow} growing perturbation (that
is, to prove theorem \ref{theorem_squares}) we need a more delicate approach.
Essentially, we cannot use approximate arithmetic progressions which start from
zero (as \( B(s,a) \) does) and the method above fails. However, translated
\( B(s,a) \)'s do exist in the neighborhood of the squares, which, as will
be shown below, is enough.

We start with an elementary estimate of Riesz products. It is interesting to
compare this to the corresponding estimate in section \ref{section_analytic_lacunary},
which requires more probabilistically oriented techniques.

\begin{lem}
\label{prod_1cos_zero}There exists a sequence \( L_{k} \) such that for every
\( \boldsymbol {\nu }=\left\{ \nu _{k}\right\} \subset \ZZ  \) satisfying
\[
\frac{\nu _{k}}{\nu _{k-1}}>L_{k}\]
 we have for almost all \( t\in \TT  \) 
\[
3^{-n}<\prod _{k=1}^{n}(1-\cos \nu _{k}t)<c^{n}\quad \forall n>N(\boldsymbol {\nu };t)\]
 where \( c<1 \) is an absolute constant.
\end{lem}
\begin{proof}
Taking \( \log  \) the problem turns into proving that
\[
-nC_{1}<\sum _{k=1}^{n}\log (1-\cos \nu _{k}t)<-nC_{2}\]
Denote 
\begin{eqnarray*}
 &  & A:=\int _{\TT }\log (1-\cos t)\, dm(t)\\
 &  & F:=\log (1-\cos t)-A\\
 &  & F_{k}:=F_{[\nu _{k}]}\quad .
\end{eqnarray*}
 One can check that for \( \{L_{k}\} \) growing sufficiently fast, 
\begin{equation}
\label{almost_orthogonal}
\left| \int F_{k}F_{k'}\right| <2^{-(k+k')}\quad \forall k\neq k'.
\end{equation}
Thus in this case, the \( F_{k} \)'s are almost orthogonal \emph{}in a sense\emph{.}
We will use the classic ``law of large numbers'' for orthonormal sequences
(see, for example, \cite{KS}, theorem 5.3.5): if \( \{\varphi _{k}\} \) is
an orthonormal system then 
\[
\sum _{k=1}^{n}\varphi _{k}(t)=o(n)\quad \mathrm{a}.\mathrm{e}.\]
The result is true also for our case. Indeed, (\ref{almost_orthogonal}) implies
the estimate
\[
\left\Vert \sum c_{k}F_{k}\right\Vert _{2}\leq C\sqrt{\sum c_{k}^{2}}\quad \forall \{c_{k}\}\]
which in turn ensures, using Schur's theorem (see \cite{KS}, theorem 3.2.2)
that the system \( F_{k} \) can be extended to an orthogonal equinormed system
on a longer interval. Therefore 
\[
\sum _{k=1}^{n}\log (1-\cos \nu _{k}t)=nA+o(n)\]
almost everywhere. The actual value of \( A \) can be calculated using the
average property for harmonic functions once it is noted that \( 1-\cos t=\frac{1}{2}|1-e^{it}|^{2} \)
and the result is 
\[
A=-\log 2\]
which proves the lemma.
\end{proof}
We may now continue with the proof of theorem \ref{theorem_squares}. For \( s,a,\nu \in \NN  \),
\( \nu >\max B(s,a) \), define 
\begin{equation}
\label{def_Bsab}
B(s,a,\nu )=\left( -\nu +B(s,a)\right) \cup \left( \nu +B(s,a)\right) \quad .
\end{equation}
Note that \( B(s,a,\nu ) \) is symmetric, because \( B(s,a) \) is.

\begin{lem}
\label{lemma_B_is_Menshov}If \( \Lambda \subset \ZZ  \) satisfies that for
some sequences \( s_{k}\rightarrow \infty ,a_{k},\nu _{k}\rightarrow \infty  \)
\begin{eqnarray*}
B_{k}:=B(s_{k},a_{k},\nu _{k}) & \subset  & \Lambda \\
\dist (B_{k},0) & \rightarrow  & \infty 
\end{eqnarray*}
then it is a Menshov spectrum.
\end{lem}
\begin{proof}
Let \( f\in L^{0} \) be any function. As before, we shall construct by induction
a series of polynomials such that the sum \( \sum P_{n} \) will converge to
\( f \). Assume \( P_{1},\ldots ,P_{n-1} \) have been constructed. Define
\begin{equation}
\label{def_fn_sparse}
F_{n}=f-\sum _{k=1}^{n-1}P_{k}\quad .
\end{equation}
We write \( \epsilon _{n}:=n^{-2} \), \( \delta _{n}:=4^{-n} \) and define
\( S_{n}=S(F_{n},n^{-2},4^{-n}) \) to be the constant defined in lemma \ref{Lemma_P_in_Bna}.
Find some \( s_{n}>S_{n} \), \( a_{n} \) and 
\begin{equation}
\label{def_nu_n}
\nu _{n}>\nu _{n-1}L_{n}
\end{equation}
 (\( L_{n} \) from lemma \ref{prod_1cos_zero}) such that 
\[
\dist (B(s_{n},a_{n},\nu _{n}),0)>\deg P_{n-1}\]
\begin{equation}
\label{def_B_sparse}
B(s_{n},a_{n},\nu _{n})\subset \Lambda 
\end{equation}
 and let \( P^{1}_{n} \) be the polynomial given by lemma \ref{Lemma_P_in_Bna}
for the same \( \epsilon _{n} \), \( \delta _{n} \) and \( F_{n} \), and
for \( s_{n} \) and \( a_{n} \). Finally, we define 
\begin{equation}
\label{def_Pn_sparse}
P_{n}(t):=\cos \nu _{n}t\cdot P^{1}_{n}(t)\quad .
\end{equation}
So, the above process is an inductive definition of sequences \( F_{n} \),
\( s_{n} \), \( a_{n} \), \( \nu _{n} \), \( P_{n}^{1} \) and \( P_{n} \)
satisfying (\ref{def_fn_sparse}), (\ref{def_nu_n}), (\ref{def_B_sparse}),
(\ref{def_Pn_sparse}) and so on. Let us now draw conclusions from these. First,
\begin{eqnarray}
 &  & \spec P_{n}\subset B(s_{n},a_{n},\nu _{n})\nonumber \\
 &  & \mathbf{m}\left\{ t\, :\, \left| F_{n}(t)\cos \nu _{n}t-P_{n}(t)\right| >4^{-n}\right\} <n^{-2}\label{Fn_cos_is_Pn} 
\end{eqnarray}
and (\ref{Fn_cos_is_Pn}) can be rephrased as
\[
\mathbf{m}\left\{ x\, :\, \left| F_{n+1}(t)-F_{n}(t)(1-\cos \nu _{n}t)\right| >4^{-n}\right\} <n^{-2}\quad .\]
So we may conclude that
\[
\left| F_{n+1}(t)-F_{n}(t)(1-\cos \nu _{n}t)\right| \leq 4^{-n}\quad \forall n>N_{1}(t).\]
Here, and in what follows, we denote by \( N_{j}(t) \) integer valued, measurable,
almost everywhere defined functions. Since \( \nu _{n}>L_{n}\nu _{n-1} \),
lemma \ref{prod_1cos_zero} implies that for almost all \( t \), 
\begin{equation}
\label{prod_cos_exponential}
3^{-n}<\prod _{k=1}^{n}|1-\cos \nu _{k}t|<c^{n}\quad \forall n>N_{2}(t)
\end{equation}
(assume that \( N_{2}>N_{1} \)). This implies a rough estimate of partial products
\[
\prod _{k=m+1}^{n}|1-\cos \nu _{k}t|<\frac{c^{n}}{3^{-m}}\quad \forall n>m>N_{2}(t)\]
 which is enough to do away with the \( 4^{-n} \) errors: if \( n>m>N_{2}(t) \)
then 
\begin{eqnarray}
\left| F_{n+1}(t)-F_{m}(t)\prod ^{n}_{k=m}(1-\cos \nu _{k}t)\right|  & \leq  & \sum _{k=m}^{n}4^{-k}\prod _{l=k+1}^{n}|1-\cos \nu _{l}t|\nonumber \\
 & < & \sum 4^{-k}\frac{c^{n}}{3^{-k}}<4c^{n}\label{Kill_errors} 
\end{eqnarray}
And therefore (using (\ref{prod_cos_exponential}))
\begin{equation}
\label{Fn_drops_expo}
|F_{n}(t)|<c_{1}^{n}\quad c_{1}<1\quad \forall n>N_{3}(t)
\end{equation}
This shows that 
\[
\sum _{k=0}^{n}P_{k}\rightarrow f\]
 so we are only left with the estimates of \( S^{*}(P_{n}) \). A little consideration
shows
\[
S^{*}(P_{n};t)\leq S^{**}(P^{1}_{n};t)\]
and condition \ref{ref_SP_lem_Bsa} of lemma \ref{Lemma_P_in_Bna} gives 
\[
\mathbf{m}\left\{ t\, :\, S^{**}(P^{1}_{n};t)>Cn^{2}(|F_{n}(t)|+4^{-n})\right\} \leq Cn^{-2}\]
so due to (\ref{Fn_drops_expo}), 
\[
S^{*}(P_{n};t)<Cn^{2}c_{1}^{n}\quad \forall n>N_{4}(t)\]
and the lemma is proved.
\end{proof}
To deduce theorem \ref{theorem_squares} consider now \( B=B(s,2a,a^{2}) \).
For sufficiently large \( a \) we will have, using proposition \ref{Bsa_linear}:
\[
b\in B\cap \ZZ ^{+}\Rightarrow b=k^{2}+\tau (b),\quad k=a+l(b)\]
with
\[
|\tau (b)|<C(s)<\sqrt{\omega (k)}\]
Lemma \ref{lemma_B_is_Menshov} implies that the set 
\[
\bigcup _{s\geq 1}B(s,2a(s),a^{2}(s))\]
is a Menshov spectrum, and the result follows.\qed

\begin{rem*}
Another direct consequence of lemma \ref{lemma_B_is_Menshov}: any symmetric
\( \Lambda  \) containing arbitrary long segments of integers is a Menshov
spectrum; first proved, using a different technique, in \cite{ayan}.
\end{rem*}

\subsection{{}}

Theorem \ref{theorem_squares} too is sharp, in the following sense:

\begin{prop}
\label{prop_squares_exact}Any sequence 
\[
\Lambda (n)=\sign (n)\cdot n^{2}+O(1)\]
 is not a Menshov spectrum.
\end{prop}
The proof will follow from the next two propositions:

\begin{prop}
\label{prop_transform_menshov}If \( \Lambda  \) is a Menshov spectrum, \( m\in \NN  \)
and \( n\in \NN  \) then
\begin{enumerate}
\item \label{ref_sub_prop_transf}\( \Lambda -n \) is a Menshov spectrum
\item \label{ref_div_prop_transf}\( \left( \frac{1}{m}\Lambda \right) \cap \ZZ  \)
is a Menshov spectrum
\end{enumerate}
\end{prop}
\begin{proof}
Start with \ref{ref_sub_prop_transf}. Let \( f\in L^{0} \). Define
\[
g(t)=e^{int}f(t)\]
 and use the fact that \( \Lambda  \) is a Menshov spectrum to write
\[
\sum _{k\in \Lambda }c_{k}e^{ikt}=g(t)\quad .\]
Clearly, this implies that \( c_{k}\rightarrow 0 \) when \( |k|\rightarrow \infty  \)
and therefore we may shift the sum in a finite number of places:
\[
\lim _{N\rightarrow \infty }\sum _{\substack {k\in \Lambda \brk |k|<N}}c_{k}e^{ikt}=\lim _{N\rightarrow \infty }\sum _{\substack {k\in \Lambda \brk |k-n|<N}}c_{k}e^{ikt}=e^{int}\cdot \lim _{N\rightarrow \infty }\sum _{\substack {l\in \Lambda -n\brk |l|<N}}c_{l+n}e^{ilt}\]
 and thus part \ref{ref_sub_prop_transf} is completed. For part \ref{ref_div_prop_transf},
define 
\[
g=f_{[m]}\]
and let \( c_{k} \) be numbers such that
\[
\sum _{k\in \Lambda }c_{k}e^{ikt}=g(t)\]
for a.e. \( t \). Then for a.e. \( t\in [0,\frac{1}{m}) \),
\[
\sum _{r=0}^{m}g(t+r/m)=\sum _{k\in \Lambda }c_{k}e^{ikt}\sum _{r}e^{ikr/m}=\sum _{\substack {k\in \Lambda \brk k\equiv 0\, (m)}}mc_{k}e^{ikt}\]
which gives
\[
f(t)=\frac{1}{m}\sum g\left( \frac{t+r}{m}\right) =\sum _{k\in \Lambda \cap m\ZZ }c_{k}e^{ikt/m}\]
and so defining \( d_{k}=c_{mk} \) we get 
\[
\sum _{k\in \left( \frac{1}{m}\Lambda \right) \cap \ZZ }d_{k}e^{ikt}=f(t)\]
 for almost every \( t\in \TT  \).
\end{proof}
In light of this, to prove proposition \ref{prop_squares_exact}, it is enough
to show that quadratic residues admit arbitrarily large holes in \( \ZZ _{p} \).
More precisely, we use the following well known fact:

\begin{prop}
\label{prop_davenport}For any \( r\geq 2 \) there exists a \( p_{0} \) such
that for any prime \( p>p_{0} \) there exists an integer \( 0\leq x<p \) satisfying
\[
\left( \frac{x+i}{p}\right) =-1\quad \forall 1\leq i\leq r\]
and \( \left( \frac{a}{p}\right)  \) is the Legendre symbol.
\end{prop}
This is a corollary of theorem 1.4.2.1 from \cite{NK}. The author attributes
this theorem to Davenport and A. Weil.

\emph{Proof of proposition \ref{prop_squares_exact}:} Let \( \Lambda  \) satisfy
\emph{}
\[
|\Lambda (n)-\sign n\cdot n^{2}|<A\quad .\]
Use proposition \ref{prop_davenport} with \( r=2A \) to get a corresponding
\( p_{0} \); pick some \( p>p_{0} \) for which \( \left( \frac{-1}{p}\right) =1 \)
(quadratic reciprocity says that this is equivalent to \( p\equiv 1 \) mod
\( 4 \)) and we shall get an interval \( (m-A,m+A) \) in \( \ZZ _{p} \) which
does not contain any \( \pm n^{2}\: (\mathrm{mod}\, p) \). Thus
\[
\Lambda \cap (m+p\ZZ )=\emptyset \]
and according to proposition \ref{prop_transform_menshov}, \( \Lambda  \)
is not a Menshov spectrum.\qed

\section{\label{section_analytic}Positive Spectra: Representation in Measure}

\subsection{{}}

We have seen in the introduction that \emph{\( \ZZ ^{+} \) is not a Menshov
spectrum}. However, replacing convergence almost everywhere by convergence in
measure enables us to get an ``analytic'' representation theorem. This result
(theorem \ref{Theorem_L0}) is proved below. Actually, we prove it in a stronger
form, inspired by Talalyan's paper \cite{survey} where the following definition
was introduced:

\begin{defn}
A series \( \sum f_{j} \) converges to \( f \) asymptotically in \( L^{p}(\TT ) \)
if for any \( \epsilon >0 \) there exists a compact \( K(\epsilon ) \) in
\( \TT  \), with \( \mathbf{m}(\TT \setminus K)<\epsilon  \), and \( \sum f_{j}=f \)
in \( L^{p}(K) \). 
\end{defn}
Clearly asymptotic convergence in \( L^{p} \) for some \( p\in [1,\infty ] \)
implies convergence in measure, while asymptotic convergence in \( L^{\infty } \)
is equivalent to convergence a.e.

We can now formulate our result.

\begin{thm}
\label{Theorem_Asympt}Given \( f\in L^{0}(\TT ) \) one can define coefficients
\( \{c_{k}\}_{k=1}^{\infty } \) such that
\[
\sum _{k>0}c_{k}e^{ikt}=f(t)\]
 asymptotically in \( L^{2}(\TT ) \).
\end{thm}

\subsection{Analytic approximation.}

The following lemma is the basis of this section. It gives an analytic approximation
of \( 1 \) with the necessary control of the partial sums.

\begin{lem}
\label{Lemma_exists_P}For every \( \epsilon >0 \) there exists an analytical
polynomial 
\[
Q(t)=Q_{\epsilon }(t)=\sum _{k>0}c_{k,\epsilon }e^{ikt}\]
 and a measurable set \( E=E_{\epsilon } \) such that 
\begin{enumerate}
\item \label{ref_Qhat_lem_AnK}\( ||\widehat{Q}||_{\infty }<\epsilon  \); 
\item \label{ref_Elarge_lem_AnK}\( \mathbf{m}(\TT \setminus E)<\epsilon  \); 
\item \label{ref_Qis1_lem_AnK}\( |Q-1|<\epsilon  \) on \( E \); 
\item \label{ref_SnQ1_lem_AnK}For every \( n \), \( ||S_{n}(Q)||_{L^{2}(E)}<2 \);
\item \label{ref_SnQ2_lem_AnK}For every \( n \), \( \mathbf{m}\left\{ t\, :\, |S_{n}(Q;t)|>2\right\} <\epsilon  \).
\end{enumerate}
\end{lem}
\begin{proof}
As in lemma \ref{Lemma_Korner}, we follow essentially the construction from
\cite{ol}. Using lemma \ref{Q_eps} fix an analytical polynomial \( G \) such
that
\begin{equation}
\label{AO_R_is_1}
||G-1||_{0}<\quarter \epsilon 
\end{equation}
and an integer \( K \) 
\begin{equation}
\label{def_K}
K>\left( \frac{2}{\epsilon }||\widehat{G}||_{1}\right) ^{2}
\end{equation}
Let \( f \) be the triangle function \( \tau _{2\pi /K} \) and approximate
it by a Fourier partial sum \( F \) such that
\begin{equation}
\label{AO_F_is_f}
||\widehat{F-f}||_{1}<\frac{\epsilon }{10K||\widehat{G}||_{1}}\quad .
\end{equation}
As before, let \( T \) be the translation by \( 2\pi /K \),
\[
T(f)(t):=f(t+2\pi K^{-1})\, ,\]
 Set
\begin{eqnarray}
N_{s} & := & K(2\deg F+\deg G+2)^{s}\label{AO_def_N} 
\end{eqnarray}
 and finally, define
\begin{eqnarray*}
Q & := & \sum _{s=1}^{K}Q_{s},\\
Q_{s} & := & T^{s}(F)\cdot G_{[N_{s}]}\quad .
\end{eqnarray*}
The choice of \( \{N_{s}\} \) ensures that all \( Q_{s} \) are analytic and
that their spectra follow one another. This implies, in particular, that 
\[
||\widehat{Q}||_{\infty }=\max ||\widehat{Q_{s}}||_{\infty }\leq \max ||\widehat{T^{s}F}||_{\infty }\max ||\widehat{G_{[N_{s}]}}||_{1}=\frac{1}{K}||\widehat{G}||_{1}<\epsilon \]
so we get requirement \ref{ref_Qhat_lem_AnK}.

Now set 
\begin{eqnarray}
I_{s} & := & \left[ \frac{2\pi (s-1)}{K},\frac{2\pi (s+1)}{K}\right] ,\qquad 1\leq s\leq K\nonumber \\
U_{s} & := & \left\{ t\in I_{s}\, :\, \left| G_{[N_{s}]}(t)-1\right| <\epsilon /4\right\} \label{AO_def_Uj} \\
E & := & \TT \setminus \bigcup U_{s}\nonumber 
\end{eqnarray}
(\ref{AO_R_is_1}) and the fact that \( K|N_{s} \) (\ref{AO_def_N}) imply
that \( \mathbf{m}U_{s}<\quarter \epsilon \cdot \mathbf{m}I_{s} \), so 
\begin{equation}
\label{E_large}
\mathbf{m}\TT \setminus E<\quarter \epsilon \quad .
\end{equation}
To verify requirement \ref{ref_Qis1_lem_AnK}, set
\[
Q':=\sum _{s=1}^{K}T^{s}(f)\cdot G_{[N_{s}]}\quad .\]
For any \( t \), for at most two values of \( s \) the value of \( T^{s}(f) \)
is non zero, and their sum is 1. Meanwhile, if \( t\in E \) then the factor
\( G_{[N_{s}]} \) belongs to \( (1-\quarter \epsilon ,1+\quarter \epsilon ) \)
for both \( s \)'s as follows from (\ref{AO_def_Uj}). So
\[
|Q'(t)-1|<\half \epsilon \quad \forall t\in E.\]
Passing from \( Q' \) to \( Q \), (\ref{AO_F_is_f}) implies
\[
\left\Vert \left( T^{s}(f-F)\cdot G_{[N_{s}]}\right) ^{\wedge }\right\Vert _{1}\leq ||\widehat{f-F}||_{1}||\widehat{G}||_{1}<\frac{\epsilon }{10K}\]
and the equality
\[
Q=Q'+\sum (T^{s}f-T^{s}F)\cdot G_{[N_{s}]}\]
 gives requirement \ref{ref_Qis1_lem_AnK}. For requirements \ref{ref_SnQ1_lem_AnK}
and \ref{ref_SnQ2_lem_AnK}, fix \( n \). Because the spectra of the \( Q_{s} \)
follow each other there exists a \( j \) such that
\[
S_{n}(Q)=\sum _{s<j}Q_{s}+S_{n}(Q_{j})\]
The sum on the right hand side can be estimated on \( E \) exactly as before
and we get

\begin{equation}
\label{AO_sum_Qs_small}
\biggl \lvert \sum _{s<j}Q_{s}(t)\biggr \rvert <1+\threequarters \epsilon \quad \forall t\in E.
\end{equation}
Whereas 
\begin{eqnarray}
||S_{n}(Q_{j})||_{L^{2}(\TT )} & \leq  & ||Q_{j}||_{2}<||T^{j}(F)||_{2}||G_{[N_{j}]}||_{\infty }<||f||_{2}||\widehat{G}||_{1}\nonumber \label{junk} \\
 & < & K^{-1/2}\cdot ||\widehat{G}||_{1}<\epsilon \label{Sn_Qj_small} 
\end{eqnarray}
and (considering \( \epsilon  \) to be sufficiently small) we get \ref{ref_SnQ1_lem_AnK}.
(\ref{Sn_Qj_small}) also gives \ref{ref_SnQ2_lem_AnK}, because we have 
\[
\mathbf{m}\left\{ t\, :\, |S_{n}(Q_{j})|>\half \right\} <4\epsilon ^{2}\]
which, together with (\ref{AO_sum_Qs_small}) and (\ref{E_large}) finishes
the lemma.
\end{proof}

\subsection{\label{ssec_proof_theorem_asympt}Proof of theorem \ref{Theorem_Asympt}.}

We shall construct by induction analytic polynomials \( P_{n} \) and sets \( E_{n} \)
such that 

\begin{enumerate}
\item \label{ref_Elarge_thm_asympt}\( \mathbf{m}(\TT \setminus E_{n})<2^{-n} \);
\item \label{ref_fissum_thm_asympt}\( \left| \sum _{k=1}^{n}P_{k}-f\right| <2^{-n} \)
on \( E_{n} \);
\item \label{ref_Pfollow_thm_asympt}\( P_{n} \) follows \( P_{n-1} \);
\item \label{ref_SkPn_thm_asympt}for any \( k \), \( ||S_{k}(P_{n})||_{L^{2}(E_{n}\cap E_{n-1})}<C2^{-n} \).
\end{enumerate}
Having these we get the theorem. Indeed, requirement \ref{ref_Pfollow_thm_asympt}
allows us to define a sequence \( \{c_{k}\}_{k=1}^{\infty } \) such that 
\[
P_{n}=\sum _{k=\deg P_{n-1}+1}^{\deg P_{n}}c_{k}e^{ikt}\quad .\]
 Requirements \ref{ref_fissum_thm_asympt} and \ref{ref_SkPn_thm_asympt} show
in a standard way that the series \( \sum c_{k}e^{ikt} \) converges to \( f \)
in \( L^{2}(\bigcap _{n=N}^{\infty }E_{n}) \) for any \( N \), and combining
this with requirement \ref{ref_Elarge_thm_asympt} we get the theorem. 

Let us now assume \( P_{1} \),...,\( P_{n-1} \) have been constructed, and
show how to build \( P_{n} \). Define
\begin{equation}
\label{PT4_def_Fn}
F_{n}:=f-\sum _{k=1}^{n-1}P_{k}
\end{equation}
 First approximate \( F_{n} \) by a trigonometric polynomial \( G_{n} \) such
that 
\begin{equation}
\label{thm_F}
|G_{n}-F_{n}|<2^{-(n+1)}\quad \mathrm{on}\: D_{n},\quad \mathbf{m}(\TT \setminus D_{n})<2^{-(n+1)}.
\end{equation}
 Put
\begin{equation}
\label{def_eps_n_thm_anal}
\epsilon (n):=\frac{2^{-(n+1)}}{||\widehat{G_{n}}||_{1}+1}
\end{equation}
and using lemma \ref{Lemma_exists_P} for \( \epsilon (n) \) we get \( Q=Q_{\epsilon (n)} \).
Set
\begin{equation}
\label{PT4_def_Pn}
P_{n}:=Q_{[r_{n}]}G_{n}
\end{equation}
where
\[
r_{n}:=\deg P_{n-1}+2\deg G_{n}+1\quad .\]
This choice of \( r_{n} \) will ensure that \( P_{n} \) are analytic, and
that \ref{ref_Pfollow_thm_asympt} is fulfilled. Define \( E_{n} \) to be 
\begin{equation}
\label{PT4_def_En}
E_{n}:=\left( E_{\epsilon (n)}\right) _{[r_{n}]}\cap D_{n}
\end{equation}
 where \( E_{\epsilon (n)} \) comes from lemma \ref{Lemma_exists_P} and \( E_{[r]} \)
is defined by the equality \( \left( \one _{E}\right) _{[r]}\equiv \one _{\left( E_{[r]}\right) } \).
Clearly, \( \mathbf{m}\left( E_{[r]}\right) =\mathbf{m}E \). So, condition
\ref{ref_Elarge_lem_AnK} of lemma \ref{Lemma_exists_P} together with (\ref{thm_F})
gives requirement \ref{ref_Elarge_thm_asympt}. Condition \ref{ref_Qis1_lem_AnK}
of lemma \ref{Lemma_exists_P} says that for \( t\in E_{n} \), \( |Q_{[r_{n}]}(t)-1|<\epsilon (n) \),
so
\begin{equation}
\label{thm_G_minus_F}
|P_{n}-G_{n}|<2^{-(n+1)}\quad \mathrm{on}\: E_{n}
\end{equation}
 and summing (\ref{PT4_def_Fn}), (\ref{thm_F}) and (\ref{thm_G_minus_F})
we end up with requirement \ref{ref_fissum_thm_asympt} fulfilled.

Requirement \ref{ref_SkPn_thm_asympt}: Now we use (\ref{Sn_of_PcrossQ}) for
the polynomial (\ref{PT4_def_Pn}):
\begin{eqnarray*}
\lefteqn {\sup _{k}||S_{k}(P_{n})||_{L^{2}(E_{n}\cap E_{n-1})}\leq } &  & \\
 & \qquad  & \leq ||G_{n}||_{L^{\infty }(E_{n}\cap E_{n-1})}\cdot \max _{k}||S_{k}(Q_{[r_{n}]})||_{L^{2}(E_{n})}+||\widehat{Q}||_{\infty }\cdot ||\widehat{G}||_{1}\\
 &  & <||G_{n}||_{L^{\infty }(D_{n}\cap E_{n-1})}\cdot \max ||S_{k}(Q)||_{L^{2}(E_{\epsilon (n)})}+2^{-n}
\end{eqnarray*}
 (here (\ref{PT4_def_En}), (\ref{def_eps_n_thm_anal}) and condition \ref{ref_Qhat_lem_AnK}
of lemma \ref{Lemma_exists_P} were used). The product on the right hand side
is
\[
<2\left\{ ||F_{n}||_{L^{\infty }(E_{n-1})}+2^{-n}\right\} \]
 due to condition \ref{ref_SnQ1_lem_AnK} of lemma \ref{Lemma_exists_P} and
(\ref{thm_F}). Finally, we use for \( F_{n} \) assumption \ref{ref_fissum_thm_asympt}
of the induction (in the \( n-1 \)'th step) and get condition \ref{ref_SkPn_thm_asympt}.
This finishes the proof of the theorem.\qed

\emph{Remarks} \textbf{\emph{}}

\begin{enumerate} 
\item[1.]In the theorem above one can replace asymptotic convergence in \( L^{2}(\TT ) \)
by the same convergence in \( L^{p}(\TT ) \) (\( p<\infty  \)). The appropriate
version of lemma \ref{Lemma_exists_P} can be proved identically, using the
classical result the \( ||S_{n}(F)||_{p}<C_{p}||F||_{p} \), \( 1<p<\infty  \),
and the proof of the theorem remains unchanged.

\item[2.] In Menshov type representation, it is often of interest to control
the speed of decrease of the the amplitudes \( \{|c_{k}|\} \), see e.g. \cite{ayan}.
The same is true for the radial representation theorem, see \cite{kk}. We wish
to pursue this direction with theorem \ref{Theorem_Asympt}. Certainly the condition
\( \{c_{k}\}\in l_{2} \) cannot be achieved here. But one can get quite close,
in a sense. Namely, either of the following conditions can be added to the formulation
of theorem \ref{Theorem_Asympt}. 

\begin{enumerate}
\item \label{cond_ck_in_lp}\( \{c_{k}\}\in l_{p}\quad \forall p>2 \)
\item \( \sum |c_{k}|^{2}w(k)<\infty  \) for any pre-given positive weight \( w(k)\searrow 0 \).
\end{enumerate}
The proof requires minor modifications. For example, let us illustrate them
for condition (a). The block structure of the ``special products'' of section
\ref{ssec_cross}, \( F\cdot G_{[r]} \) for \( r>2\deg F \), implies that
\begin{equation}
\label{fnorm_of_cross}
||\widehat{F\cdot G_{[r]}}||_{p}=||\widehat{F}||_{p}\cdot ||\widehat{G}||_{p}
\end{equation}
Using that, a direct inspection of the proof of theorem \ref{Theorem_Asympt}
shows that it is enough, given \( p>2 \) to replace condition \ref{ref_Qhat_lem_AnK}
in lemma \ref{Lemma_exists_P} with the stronger\\
\ref{ref_Qhat_lem_AnK}\( ' \) \( ||\widehat{Q}||_{p}<\epsilon  \).\\
This, in turn, comes from the following estimate for the triangle function
\[
||\widehat{\tau _{h}}||^{p}_{p}=O(h^{p-1})\quad .\]
Which gives
\[
||\widehat{Q_{s}}||_{p}=||\widehat{T^{s}(F)}||_{p}\cdot ||\widehat{G_{[N_{s}]}}||_{p}\leq ||\widehat{\tau _{2\pi /K}}||_{p}\cdot ||\widehat{G}||_{p}=O\left( K^{\frac{1-p}{p}}\right) \cdot C(\epsilon )\]
and the separation of spectra gives
\[
||\widehat{Q}||_{p}^{p}=\sum _{s=1}^{K}||\widehat{Q_{s}}||_{p}^{p}=O(K^{-p})\cdot C(\epsilon )\]
and replacing (\ref{def_K}) with a suitable sufficiently small \( K \), we
get \ref{ref_Qhat_lem_AnK}\( ' \).

\end{enumerate}

\subsection{Representation of infinity.}

In \cite{12}, Menshov proved that for convergence in measure one can extend
the class of ``representable functions'' to include ones taking infinite values.
It is natural to formulate this result for real functions. Let \( f \) be a
measurable function from \( \TT  \) into \( \RR \cup \{\pm \infty \} \). Menshov's
result says that any such \( f \) can be expanded to a real trigonometric series:
\begin{equation}
\label{real_series}
f(t)=\sum _{k\in \ZZ }d(k)e^{ikt}\quad d(-k)=\overline{d(k)}
\end{equation}
converging in measure. Theorem \ref{Theorem_L0} admits also functions as above,
and thus gives a direct improvement of Menshov's result. More precisely the
following proposition is true: 
\theoremstyle{plain} 
\newtheorem*{thm3p}{Theorem \ref{Theorem_L0}$'$}
\begin{thm3p}
 Any \( f \) as above can be represented by 
\begin{equation}
\label{S44_expansion}
f(t)=\sum _{k>0}c(k)e^{ikt}
\end{equation}
where the series converges in measure. 
\end{thm3p} Convergence here means that the real and imaginary parts of the
series converge respectively to \( f \) and \( 0 \) in measure.

We first mention that having (\ref{S44_expansion}) one immediately gets the
``real'' expansion (\ref{real_series}) by taking \( d(k) \) to be \( \half c(k) \)
for \( k>0 \) and \( \half \overline{c(k)} \) for \( k<0 \). To get (\ref{S44_expansion})
we denote by \( A \), \( B^{+} \), \( B^{-} \) the sets where \( f \) takes
finite, \( +\infty  \) and \( -\infty  \) values, respectively, and define
\[
f_{n}:=\left\{ \begin{array}{ll}
f & \mathrm{on}\: A\\
n & \mathrm{on}\: B^{+}\\
-n & \mathrm{on}\: B^{-}
\end{array}\right. .\]
Now, keeping conditions \ref{ref_Elarge_thm_asympt} and \ref{ref_Pfollow_thm_asympt}
in \ref{ssec_proof_theorem_asympt} unchanged, we replace \( f \) by \( f_{n} \)
in \ref{ref_fissum_thm_asympt} (and in (\ref{PT4_def_Fn})) and replace condition
\ref{ref_SkPn_thm_asympt} by these two: 
\begin{enumerate} 
\item[$\mbox{\ref{ref_SkPn_thm_asympt}}_1$] \( \max _{k}||S_{k}(P_{n})||_{L^{2}(E_{n}\cap E_{n-1}\cap A)}<C2^{-n} \)

\item[$\mbox{\ref{ref_SkPn_thm_asympt}}_2$] \( \max _{k}||S_{k}(P_{n})||_{L^{2}(E_{n}\cap E_{n-1}\cap (B^{+}\cup B^{-}))}<C \)

\end{enumerate} Repeating the arguments in \ref{ssec_proof_theorem_asympt}
literally one constructs polynomials \( \{P_{n}\} \) satisfying the modified
conditions. It is easy to see that these conditions imply that the series (\ref{S44_expansion})
yielded by \( \sum _{k>0}P_{k} \) converges to \( f \) in measure.

\section{\label{section_analytic_lacunary}Positive sparse spectra}

\subsection{{}}

The aim of this section is to improve the results of the previous one, showing
that only relatively ``small'' parts of \( \ZZ ^{+} \) are needed. For brevity,
we shall discuss only convergence in measure, rather than asymptotic \( L^{p} \),
and only finite functions.

\begin{defn}
A set \( \Lambda  \) of integers is a \emph{Menshov spectrum} \emph{in measure}
if every \( f\in L^{0}(\TT ) \) admits a representation 
\[
f(t)=\sum _{k\in \Lambda }c_{k}e^{ikx}\equiv \lim _{N\rightarrow \infty }\sum _{k\in \Lambda ,|k|\leq N}c_{k}e^{ikx}\]
converging in measure.
\end{defn}
We shall prove the following analogues of the results of section \ref{section_lacunary}:

\begin{thm}
\label{theorem_analytic_sparse}Given a positive sequence \( \epsilon (n)=o(1) \)
as \( n\rightarrow \infty  \), one can define a Menshov spectrum in measure
\( \Lambda \subset \ZZ ^{+} \) such that
\begin{equation}
\label{analytic_lambda_exp}
\frac{\lambda (n+1)}{\lambda (n)}>1+\epsilon (n)\quad n=1,2,\ldots 
\end{equation}
 
\end{thm}
{}

\begin{thm}
\label{theorem_analytic_arutyunyan}If a set \( \Lambda \subset \ZZ ^{+} \)
contains arbitrarily long segments then it is a Menshov spectrum in measure.
\end{thm}
{}

\begin{thm}
\label{theorem_analytic_squares}For any sequence \( w(k)\rightarrow \infty  \)
as \( k\rightarrow \infty  \) one can construct a Menshov spectrum in measure
\( \Lambda  \),
\[
\Lambda =\left\{ k^{2}+o(w(k))\right\} \quad .\]

\end{thm}

\subsection{The blocks \protect\( D(s,a)\protect \)}

The equivalent of the blocks \( B(s,a) \) from the proof of theorem \ref{Theorem_sparse}
will be denoted by \( D(s,a) \) and defined simply by
\[
D(s,a):=B^{+}_{1}(s,(2s)^{2s+2}a)+B_{2}(s,a)\]
where \( B^{+}_{1} \) and \( B_{2} \) are from the definition of \( B(s,a) \).
Of course, \( D(s,a)\subset \ZZ ^{+} \). The substitute for lemma \ref{Lemma_P_in_Bna}
is the following:

\begin{lem}
\label{lemma_analytic_bsa}For every \( \epsilon >0 \) and \( f\in L^{0} \)
there exists an \( S=S(f,\epsilon ) \) with the following property. Given \( s>S \)
and \( a\in \NN  \) one may construct a trigonometric polynomial \( P \) satisfying 
\begin{enumerate}
\item \label{req_fisP_lem_ABsa}\( ||f-P||_{0}<\epsilon  \);
\item \label{req_specP_lem_ABsa}\( \spec P\subset D(s,a) \);
\item \label{req_SnP_lem_ABsa}For any \( n \), \( \mathbf{m}\left\{ t\, :\, |S_{n}(P;t)|>2|f(t)|+\epsilon \right\} <\epsilon  \).
\end{enumerate}
\end{lem}
\begin{proof}
The proof of this lemma is similar to that of lemma \ref{Lemma_P_in_Bna}, replacing
the use of lemma \ref{Lemma_Korner} with lemma \ref{Lemma_exists_P}, so we
will describe it briefly.

First, approximate \( f \) by a polynomial \( P_{1} \) satisfying 
\begin{equation}
\label{VE_f_minus_p1}
||f-P_{1}||<\epsilon _{1}:=\third \epsilon 
\end{equation}
 then use lemma \ref{Q_eps} to get an analytic polynomial \( Q_{2} \) satisfying
\begin{equation}
\label{VE_eps_2_del_2}
||Q_{2}-1||_{0}<\epsilon _{2}:=\frac{\epsilon }{6\max \{||\widehat{P_{1}}||_{1},\deg P_{1}\}+3}\quad .
\end{equation}
 Now use lemma \ref{Lemma_exists_P} with 
\begin{equation}
\label{VE_eps_3_del_3}
\epsilon _{3}:=\frac{\epsilon }{6||\widehat{P_{1}}||_{1}||\widehat{Q_{2}}||_{1}+3}
\end{equation}
to get an analytic \( Q_{3} \) satisfying

\begin{eqnarray}
 &  & ||\widehat{Q_{3}}||_{\infty }\leq \epsilon _{3};\label{VE_Q3hat_is_small} \\
 &  & \mathbf{m}\left\{ t\, :\, |S_{n}(Q_{3};t)|>2\right\} <\epsilon _{3}\quad \forall n;\label{VE_Sn_Q3} \\
 &  & ||Q_{3}||_{0}<\epsilon _{3}\quad .\nonumber 
\end{eqnarray}
 Set 
\[
S(f,\epsilon ):=\max \left\{ \deg P_{1},\deg Q_{2},\deg Q_{3}\right\} \quad .\]
Fix some \( s>S \) and any integer \( a \). Defining, as in (\ref{def_P2})
\begin{eqnarray}
P_{2} & := & \sum _{|k|\leq \deg P_{1}}\widehat{P_{1}}(k)e^{ikt}\cdot (Q_{2})_{[p(k)]}\label{VE_def_P2} \\
p(k) & := & a(2s)^{k+s}\nonumber 
\end{eqnarray}
 we will have \( \spec P_{2}\subset B_{2}(s,a) \). 

The same estimate as in (\ref{P1_minus_P2}) gives 
\begin{equation}
\label{VE_f_minus_P2}
||P_{2}-f||_{0}<\twothirds \epsilon \quad .
\end{equation}
 Finally, define
\[
P\equiv P_{3}=(Q_{3})_{[p(s+2)]}P_{2}\quad .\]
The fact that \( Q_{3} \) is analytic gives that \( \spec (Q_{3})_{[p(s+2)]}\subset B^{+}_{1}(s,p(s+2)) \),
which implies \ref{req_specP_lem_ABsa}. As in lemma \ref{Lemma_P_in_Bna},
we have
\[
\mathbf{m}\left\{ x\, :\, |P-P_{2}|>\third \epsilon \right\} \leq \mathbf{m}\left\{ t\, :\, |1-Q_{3}|>\frac{\epsilon }{3||P_{2}||_{\infty }}\right\} \quad ,\]
\begin{equation}
\label{VE_P2_not_too_large}
||P_{2}||_{\infty }\leq ||\widehat{P_{2}}||_{1}\leq ||\widehat{P_{1}}||_{1}||\widehat{Q_{2}}||_{1}
\end{equation}
so
\begin{equation}
\label{VE_P3_minus_P2}
\mathbf{m}\left\{ t\, :\, |P-P_{2}|>\third \epsilon \right\} \leq \mathbf{m}\left\{ t\, :\, |1-Q_{3}|>\epsilon _{3}\right\} <\third \epsilon \quad .
\end{equation}
Summing (\ref{VE_f_minus_P2}) and (\ref{VE_P3_minus_P2}) gives \ref{req_fisP_lem_ABsa}.
To estimate \( S_{n}(P) \) use (\ref{p_s_plus_2_large}) to write, as in section
\ref{ssec_cross}, \( n=r\cdot p(s+2)+l \) and get
\[
S_{n}(P;t)\leq |S_{r}(Q_{3};p(s+2)t)|\cdot |P_{2}(t)|+2||\widehat{Q_{3}}||_{\infty }||\widehat{P_{2}}||_{1}\]
as before, estimate \( ||\widehat{P_{2}}||_{1} \) by (\ref{VE_P2_not_too_large}).
Now (\ref{VE_Sn_Q3}) and (\ref{VE_eps_3_del_3}) give
\[
S_{n}(P;t)\leq 2|P_{2}(t)|+\third \epsilon \quad .\]
outside a set of measure \( \epsilon _{3}<\frac{1}{3}\epsilon  \), and remembering
(\ref{VE_f_minus_P2}), the lemma is proved.
\end{proof}
Lemma \ref{lemma_analytic_bsa} proves theorem \ref{theorem_analytic_sparse}
in much the same way as for theorem \ref{Theorem_sparse}.

\subsection{``Analytic'' Riesz products.}

Passing to the proofs of theorems \ref{theorem_analytic_arutyunyan} and \ref{theorem_analytic_squares}
we make the following remark. In the ``squares'' theorem of section \ref{section_lacunary},
given a trigonometric polynomial, we used multiplication with a Riesz product
to localize the spectrum in a corresponding block of high frequencies and the
estimate (\ref{prod_cos_exponential}) provided the necessary control of partial
sums. In turn, that estimate came essentially from the inequality
\[
\int _{\TT }\log (1-\cos t)\, dm<0\quad .\]
 Now only the positive part of the spectrum is available, so Riesz products
should be replaced by its ``analytic'' counterpart
\begin{equation}
\label{def_qn}
q_{n}(t):=\prod _{k\leq n}(1-e^{i\nu (k)t})\quad .
\end{equation}
However, 
\[
\int _{\TT }\log |1-e^{it}|=0\]
and the sequence \( \{|q_{n}|\} \) indeed does not tend to zero. But it repeatedly
returns close to zero for almost every \( t \). Fortunately, this weaker property
is enough to manage the localization of spectrum with a simultaneous control
of the partial sums, corresponding to convergence in measure.

\begin{lem}
\label{lemma_exp_is_martingale}There exists numbers \( L_{k} \) such that
for every \( \{\nu _{k}\}\subset \ZZ  \) satisfying
\begin{equation}
\label{G2_def_ni}
\frac{\nu _{k}}{\nu _{k-1}}>L_{k}
\end{equation}
 one has 
\begin{equation}
\label{liminf_prod}
\liminf _{N\rightarrow \infty }\prod _{k=1}^{N}|1-e^{i\nu _{k}t}|=0\quad \mathrm{a}.\mathrm{e}.
\end{equation}

\end{lem}
\begin{proof}
This is equivalent to the following quality
\begin{equation}
\label{liminf_sum}
\liminf _{N\rightarrow \infty }\sum _{k=1}^{N}\log |1-e^{i\nu _{k}t}|=-\infty \quad \mathrm{a}.\mathrm{e}.
\end{equation}
 This is a direct consequence of a well known principle that sparse subsequences
of a system \( \left\{ \varphi (\nu t)\right\}  \) inherit properties of independent
variables, and in our case, that the sum
\[
\sum _{k=1}^{N}\varphi (\nu _{k}t)\]
behaves like a random walk for a sufficiently fast growing \( \nu _{k} \).
We use the following version of this fact, which can be found (in a more general
form) in the paper of Gaposhkin \cite{15} (corollary 1): if \( \int _{\TT }\varphi \, dm=0 \),
\( ||\varphi ||_{L^{2}(\TT )}=1 \) and
\begin{equation}
\label{G2_omega}
\sum _{j<k}\omega _{2}\left( \frac{n_{j}}{n_{k}};\varphi \right) =o(1),\qquad \sum _{j>k}\omega _{2}\left( \frac{n_{k}}{n_{j}};\varphi \right) =o(1)
\end{equation}
where \( \omega _{2} \) is the module of continuity of \( \varphi  \) in the
space \( L^{2} \), then for every set \( E\subset \TT  \), \( mE>0 \), 
\begin{equation}
\label{G2_clt}
\mathbf{m}\left\{ t\in E\, :\, \frac{1}{\sqrt{N}}\sum _{k=1}^{N}\varphi (\nu _{k}t)<y\right\} \rightarrow \mathbf{m}E\cdot \frac{1}{\sqrt{2\pi }}\int _{-\infty }^{y}e^{-t^{2}/2}dt,\quad N\rightarrow \infty 
\end{equation}
Taking \( \varphi  \) to be \( \log |1-e^{it}| \) (after a suitable normalization)
one can see that for \( L_{k} \) growing sufficiently fast, (\ref{G2_def_ni})
implies (\ref{G2_omega}). If (\ref{liminf_sum}) fails on a set \( E \) of
positive measure we get an immediate contradiction to (\ref{G2_clt}) e.g. for
\( y=0 \).
\end{proof}
\begin{rem*}
\label{remark_low_anal_est}As before, we shall need a rough estimate from below,
\begin{equation}
\label{qn_larger_exp}
\left( \threequarters \right) ^{n}<\prod _{k=1}^{n}|1-e^{i\nu _{k}t}|\quad \forall n>K_{1}(t)
\end{equation}
where, as in the proof of theorem \ref{theorem_squares}, here and below \( K_{j}(t) \)
are measurable functions \( \TT \rightarrow \NN  \) defined almost everywhere.
This estimate can be achieved exactly as in lemma \ref{prod_1cos_zero}.
\end{rem*}

\subsection{The main lemma.}

In this section, we shall combine (\ref{liminf_prod}) with stopping time techniques,
to prove an analogue of lemma \ref{Bna_is_menshov}, lemma \ref{lemma_D_is_Menshov}
below. Let us start with a definition:

\begin{defn}
For \( s \), \( a \), \( \nu \in \ZZ ^{+} \) we define 
\[
D(s,a,\nu )=\nu +D(s,a)\]
and call a set \( \Lambda  \) a ``\( D \)-set'' if one has a sequence \( D(s_{k},a_{k},\nu _{k})\subset \Lambda  \)
with \( s_{k} \) and \( \nu _{k} \) tending to \( \infty  \).
\end{defn}
\begin{lem}
\label{lemma_AO_stoptime}Let \( \Lambda \subset \ZZ ^{+} \) be a \( D \)-set;
\( f\in L^{0}(\TT ) \) be any function vanishing outside of an interval \( I \);
and \( \epsilon >0 \) be some number. Then there exists a polynomial \( P \)
with the following properties:
\begin{enumerate}
\item \label{req_fisP_lem_stoptime}\( ||P-f||_{0}<\epsilon  \);
\item \label{req_specP_lem_stoptime}\( \spec P\subset \Lambda  \) ;
\item \label{req_SnP_lem_stoptime}For any \( n\in \ZZ ^{+} \), 
\[
\mathbf{m}\left\{ t\notin I\, :\, |S_{n}(P;t)|>\epsilon \right\} <\epsilon \qquad .\]

\end{enumerate}
\end{lem}
\begin{proof}
1. We shall construct inductively a sequence of polynomials \( P_{k} \) using
the following procedure:
\begin{eqnarray}
T_{k} & := & f-\sum _{l=0}^{k-1}P_{l}\label{Def_Tk} \\
R_{k}(t) & := & \left\{ \begin{array}{ll}
T_{k}(t) & |T_{l}(t)|>\half \epsilon \: \forall l\leq k\: \mathrm{and}\: t\in I\\
0 & \mathrm{otherwise}
\end{array}\right. \quad .\label{Rk_pp_is_stop_time} 
\end{eqnarray}
We now use lemma \ref{lemma_analytic_bsa} for the function \( R_{k} \) and
the parameter \( \epsilon 2^{-k-2} \). This gives us some number \( S=S(R_{k},\epsilon 2^{-k-2}) \).
let \( L_{k} \) be numbers from lemma \ref{lemma_exp_is_martingale} and choose
\( a_{k} \), \( \nu _{k} \) and \( s_{k} \) to satisfy
\begin{eqnarray}
D(s_{k},a_{k},\nu _{k}) & \subset  & \Lambda \label{sn_put_Pn_in_Lambda} \\
\nu _{k} & > & L_{k}\nu _{k-1}\label{exp_sk_martingale} \\
s_{k} & > & S\quad .\label{sk_larger_S} 
\end{eqnarray}
Returning to lemma \ref{lemma_analytic_bsa}, (\ref{sk_larger_S}) allows us
to define a polynomial \( P_{k}' \) satisfying: 
\begin{eqnarray}
\spec P_{k}' & \subset  & D(s_{k},a_{k})\label{Pk_p_in_Dsa} \\
||R_{k}-P_{k}'||_{0} & < & \epsilon 2^{-k-2}\label{Pk_p_is_Rk_p} 
\end{eqnarray}
 and for any \( n \), 
\begin{equation}
\label{Sn_Pk_p_small}
\mathbf{m}\left\{ t\, :\, |S_{n}(P_{k}';t)|>2|R_{k}(t)|+\epsilon 2^{-k-2}\right\} <\epsilon 2^{-k-2}
\end{equation}
 Finally, we define 
\begin{equation}
\label{Pk_is_with_exp}
P_{k}:=e^{i\nu _{k}t}\cdot P_{k}'\quad .
\end{equation}
So, the above process is an inductive definition of sequences \( T_{k} \),
\( R_{k} \), \( s_{k} \), \( a_{k} \), \( \nu _{k} \), \( P_{k}' \) and
\( P_{k} \) satisfying (\ref{Def_Tk})-(\ref{Pk_is_with_exp}). 

2. First we claim 
\begin{equation}
\label{Tk_to_zero}
||T_{k}||_{0}<\epsilon \quad \forall k>K(\epsilon )
\end{equation}
 To see this, collect (\ref{Rk_pp_is_stop_time}), (\ref{Pk_p_is_Rk_p}) and
(\ref{Pk_is_with_exp}) and get the following for \( T_{k} \):
\begin{equation}
\label{Rk_p_recursion}
T_{k+1}(t)=\left\{ \begin{array}{ll}
T_{k}(t)(1-e^{i\nu _{k}t})+r_{k}(t) & |T_{l}(t)|>\half \epsilon \: \forall l\leq k\: \mathrm{and}\: t\in I\\
T_{k}(t)+r_{k}(t) & \mathrm{otherwise}
\end{array}\right. 
\end{equation}
with 
\begin{equation}
\label{Ek_is_small}
||r_{k}||_{0}<\epsilon 2^{-k-2}
\end{equation}
In particular, 
\begin{equation}
\label{Ek_small}
|r_{k}(t)|<2^{-k-2}\quad \forall k>K_{2}(t)\quad \mathrm{a}.\mathrm{e}.
\end{equation}
 Denoting 
\[
q_{n}:=\prod _{k=1}^{n}(1-e^{i\nu _{k}t})\quad ,\]
 we have (\ref{qn_larger_exp}) and (\ref{liminf_prod}) (remember (\ref{exp_sk_martingale}))
i.e.
\begin{equation}
\label{prod_exp_high}
|q_{n}|>\left( \threequarters \right) ^{n}\quad \forall n>K_{1}(t);
\end{equation}
\begin{equation}
\label{prod_exp_low}
\liminf _{n\rightarrow \infty }|q_{n}(t)|=0\quad \mathrm{a}.\mathrm{e}.
\end{equation}
 Set:
\begin{equation}
\label{def_V_eps}
V_{\epsilon }:=\{t\in I\, :\, |T_{k}(t)>\half \epsilon \quad \forall k\in \NN \}
\end{equation}
i.e. the set where \( T_{k} \) never ``stabilizes''. Using (\ref{Rk_p_recursion}),
(\ref{Ek_small}) and (\ref{prod_exp_low}) we have for almost every \( t\in V_{\epsilon } \)
and \( k>k_{1}>\max \{K_{1}(t),K_{2}(t)\} \),
\begin{eqnarray}
\left| T_{k+1}(t)-T_{k_{1}}(t)\prod ^{k}_{l=k_{1}}(1-e^{i\nu _{l}t})\right|  & \leq  & \sum _{l=k_{1}}^{k}2^{-l-2}\prod _{j=l+1}^{k}|1-e^{i\nu _{j}t}|\nonumber \\
 & < & \sum _{l\leq k}2^{-l-2}\left( \threequarters \right) ^{l}\cdot |q_{k}(t)|\nonumber \label{how_do_i_get_rid_of_this_2} 
\end{eqnarray}
so
\[
\left| T_{k+1}(t)-\frac{T_{k_{1}}(t)}{q_{k_{1}-1}(t)}q_{k}(t)\right| <C|q_{k}(t)|\]
contrasting this with (\ref{prod_exp_high}) and (\ref{def_V_eps}) we conclude
\begin{equation}
\label{m_V_eps_zero}
\mathbf{m}V_{\epsilon }=0\quad .
\end{equation}
Now, (\ref{m_V_eps_zero}) gives that for some \( K_{3}(t) \),
\[
T_{k+1}(t)=T_{k}(t)+r_{k}(t)\quad \forall k\geq K_{3}(t)\]
 and
\[
|T_{K_{4}(t)}|\leq \half \epsilon \]
Choosing \( K(\epsilon ) \) sufficiently large to have
\[
\mathbf{m}\left\{ t\, :\, K_{4}(t)>K(\epsilon )\right\} <\half \epsilon \]
and using (\ref{Ek_is_small}), we get (\ref{Tk_to_zero}). 

3. Let
\[
P:=\sum _{k=0}^{K(\epsilon )}P_{k}\quad .\]
(\ref{Tk_to_zero}) and (\ref{Def_Tk}) give \ref{req_fisP_lem_stoptime}. (\ref{sn_put_Pn_in_Lambda}),
(\ref{Pk_p_in_Dsa}), and (\ref{Pk_is_with_exp}) give \ref{req_specP_lem_stoptime}.
To see \ref{req_SnP_lem_stoptime}, note that 
\[
S_{n}(P_{k})=\left\{ \begin{array}{ll}
e^{i\nu _{k}t}\cdot S_{n-\nu _{k}}(P_{k}') & \forall n\geq \nu _{k}\\
0 & \mathrm{otherwise},
\end{array}\right. \]
 (\ref{Sn_Pk_p_small}) means that 
\begin{equation}
\label{Sn_Pk_small}
\mathbf{m}\left\{ t\, :\, |S_{n}(P_{k};t)|>2|R_{k}(t)|+\epsilon 2^{-k-2}\right\} <\epsilon 2^{-k-2}\quad .
\end{equation}
Since \( t\notin I\Rightarrow R_{k}(t)=0 \), summing up the last inequality
over \( k \) we get \ref{req_SnP_lem_stoptime}.
\end{proof}
\begin{lem}
\label{lemma_D_is_Menshov}Any \( D \)-set is a Menshov spectrum in measure.
\end{lem}
\begin{proof}
Let \( \Lambda \subset \ZZ ^{+} \) be a \( D \)-set. Let \( f\in L^{0}(\TT ) \)
be any function. Take a sequence of intervals \( \{I_{k}\} \) with 
\[
|I_{k}|\rightarrow 0\]
which cover each point of the circle infinitely many times. We shall construct
inductively a sequence of polynomials \( P_{k} \) using the following procedure.
Define
\[
R_{k}:=\one _{I_{k}}\cdot \left( f-\sum _{l=0}^{k-1}P_{l}\right) \]
and use lemma \ref{lemma_AO_stoptime} with the data 
\[
\Lambda _{k}:=\Lambda \cap [\deg P_{k-1},\infty )\]
(which is clearly a \( D \)-set), \( R_{k} \), \( I_{k} \) and \( \epsilon =2^{-k} \).
We get a polynomial \( P_{k} \) which satisfies:
\begin{eqnarray}
\spec P_{k} & \subset  & \Lambda _{k}\subset \Lambda \nonumber \\
||P_{k}-R_{k}||_{0} & < & 2^{-k}\label{Pk_is_Rk_lemmaD} 
\end{eqnarray}
 and for any \( n\in \NN  \), 
\begin{equation}
\label{Sn_Pk_small_lemmaD}
\mathbf{m}\left\{ t\notin I_{k}\, :\, |S_{n}(P_{k};t)|>2^{-k}\right\} <2^{-k}\qquad .
\end{equation}
 The sum of the \( P_{k} \)'s, written as a trigonometric series will prove
the lemma. As it is clearly supported on \( \Lambda  \), we need only show
that it converges to \( f \) in measure. One can easily see that (\ref{Pk_is_Rk_lemmaD})
and the fact that the \( I_{k} \)'s cover \( \TT  \) infinitely many times
implies that 
\[
\sum _{k=0}^{\infty }P_{k}(t)=f(t)\quad \mathrm{a}.\mathrm{e}.\]
 As the \( P_{k} \)'s follow each other, we write
\[
P_{k}=\sum _{n=\deg P_{k-1}+1}^{\deg P_{k}}c_{n}e^{int}\]
and use (\ref{Sn_Pk_small_lemmaD}), written as
\[
\max _{n}||S_{n}(P_{k})||_{0}<2^{-k}+|I_{k}|\rightarrow 0\quad k\rightarrow \infty \]
to get that \( \sum _{n\in \Lambda }c_{n}e^{int}=f(t) \), in measure.
\end{proof}

\subsection{{}}

Theorem \ref{theorem_analytic_arutyunyan} follows from this lemma immediately
while theorem \ref{theorem_analytic_squares} is derived in much the same way
that theorem \ref{theorem_squares} follows from lemma \ref{lemma_B_is_Menshov}.

Finally, we mention that theorems \ref{theorem_analytic_sparse} and \ref{theorem_analytic_squares}
are sharp in the same sense that theorems \ref{Theorem_sparse} and \ref{theorem_squares}
respectively are. Theorem \ref{theorem_analytic_sparse} is sharp for exactly
the same reasons as theorem \ref{Theorem_sparse}; and since proposition \ref{prop_transform_menshov}
clearly holds (with an identical proof) for Menshov spectra in measure, we may
conclude the sharpness of theorem \ref{theorem_analytic_squares} using proposition
\ref{prop_davenport}.

\end{document}